\documentclass{article}

\usepackage{amssymb}
\usepackage{dsfont}

\newtheorem{lem}{Lemma}

\newtheorem{theo}[lem]{Theorem}

\newtheorem{fait}[lem]{Fact}
\newtheorem{cor}[lem]{Corollary}

\newcommand{\proof}{\noindent{\bf Proof.}~}
\newcommand{\qed}{\ \hfill$\square$\bigskip}

\newcommand{\rk}{\hbox{\rm rk}\,}

\newcommand{\SL}{\hbox{\rm SL}\,}
\newcommand{\PSL}{\hbox{\rm PSL}\,}

\renewcommand{\>}{\rangle}
\renewcommand{\o}{^\circ}

\title{Cosets, genericity, and the Weyl group}
\author{Eric Jaligot}
\date{\today}

\begin{document}
\maketitle

\begin{abstract}
In a connected group of finite Morley rank in which, generically, elements 
belong to connected nilpotent subgroups, proper normalizing cosets of definable 
subgroups are not generous. We explain why this is true and what consequences 
this has on an abstract theory of Weyl groups in groups of finite Morley rank.
\end{abstract}

The only known infinite simple groups of finite Morley rank are the
simple algebraic groups over algebraically closed fields and this is a motivation, 
among many others, for a classification project of these groups. 
It borrows ideas and techniques from the Classification
of the Finite Simple Groups but at the same time it may provide, sometimes, 
a kind of simplified version of the finite case. This is mostly due to the
existence of well-behaved notions of genericity and connectivity in the
infinite case, which unfortunately have no direct finite analogs.

The present note deals with a very specific and technical topic concerning such
arguments based on genericity in the case of infinite groups of finite
Morley rank, which serve here to bypass allegro potential complications of various nature, 
including finite combinatorics. 
As a result, we show similarities with algebraic groups in any case as far
as a theory of Weyl groups is concerned, and naturally this applies also 
to non-algebraic configurations which are encountered throughout much of the 
current work in the area. 

In a connected reductive algebraic group, maximal (algebraic) tori are conjugate
and cover the group generically, with the Weyl group governing essentially
the structure of the entire group. In the abstract context, we use the term 
``generous" to speak of a subset 
``whose union of conjugates is generic in the group", the typical property of tori in 
the classical algebraic case. There are at least two abstract
versions of tori in groups of finite Morley rank, which coincide at least in the case
of a reductive algebraic group, {\em decent tori} on the one hand
and {\em Carter subgroups} on the other. The main caveat with these two
more abstract notions for a seemingly complete analogy with algebraic groups 
is in both cases an unknown existence, more precisely the existence of a 
{\em nontrivial} decent torus on the one hand and the existence of 
a {\em generous} Carter subgroup on the other. Anyway, here we follow 
an approach resolutely adapted to the second notion. 

With both notions there are conjugacy theorems, the conjugacy of 
maximal decent tori \cite{Cherlin05} and of generous Carter subgroups \cite{Jaligot06}. 
This gives a natural notion of Weyl group in each case, 
$N(T)/C\o(T)$ for some maximal decent torus $T$ or $N(Q)/Q$ for some
generous Carter subgroup $Q$. In any case and whatever the Weyl group
is, it is finite and, as with classical Weyl groups and algebraic tori in
algebraic groups, its determination and its action on the underlying subgroup is
fundamental in the abstract context.

As an element of the Weyl group is a coset in the ambient
group, it is then useful to get a description of such cosets, even though recovering from 
such a description the structure and the action of the Weyl group is in general 
a particularly delicate task. This is mostly due to the fact that, in practice, one 
can only get a generic, and thus weak, description of the coset. 
In \cite{CherlinJaligot2004} such arguments were however developed intensively,
and this was highly influenced by one of the most critical aspects of the early 
work, notably by Nesin, on the so-called ``bad" groups of finite Morley rank 
(\cite[Theorem 13.3]{BorovikNesin(Book)94}).
In this paper, a pathological coset, whose representative is typically a
Weyl element which should not exist, is usually shown to be both generous {\em and}
nongenerous, and then the coset does not exist. This is the main
protocol, sometimes refered to as ``coset arguments", for the
limitation of the size of the Weyl group. Generosity is usually obtained
by unexpected commutations between the Weyl elements and the underlying
subgroup, and in general this may depend on the specific configuration considered. It is
certainly the pathological property in any case, and we shall prove here
at a reasonable level of generality that the existing cosets should be nongenerous. 

In particular, we rearrange as follows the protocol of \cite{CherlinJaligot2004} 
in the light of further developments of \cite{Jaligot06} concerning generosity.

\begin{theo}[Generix and the Cosets]\label{TheoGenerixCosets}
Let $G$ be a connected group of finite Morley rank in which, generically, elements 
belong to connected nilpotent subgroups. Then the coset $wH$ is not generous for any 
definable subgroup $H$ and any element $w$ normalizing $H$ but not in $H$.
\end{theo}

The assumption on the generic elements of $G$ in Theorem \ref{TheoGenerixCosets} 
can take several forms, and we will explain this shortly. 
The most typical case where Theorem \ref{TheoGenerixCosets} applies 
is however the case in which $H=Q$ is a generous Carter subgroup of $G$. 
In particular, the present paper is also an appendix of \cite{Jaligot06} on the 
structure of groups of finite Morley rank with such a generous Carter subgroup, 
and more precisely a follow-up to Section 3.3 in that paper.

The general idea of the protocol of \cite{CherlinJaligot2004} has been used repeatedly in
various contexts, most notably to get a fine description of
$p$-torsion in terms of connected nilpotent subgroups of bounded exponent and of decent
tori \cite{BurdgesCherlinSemisimpleTorsion}. 
Applied to the most natural kind of Weyl groups, the protocol shows that 
centralizers of decent tori are connected in any connected group, 
implying in particular that the Weyl group $N(T)/C\o(T)$ attached to a decent 
torus $T$ acts faithfully on $T$. This corresponds to the most typical and smooth 
applications of the protocol in \cite{CherlinJaligot2004}, generally a lemma expedited at the 
early stage of the analysis of each configuration considered there. With 
\cite{Cherlin05} and \cite{Jaligot06}, and eventually the 
finiteness of conjugacy classes of uniformly definable families of decent tori 
of \cite[Theorem 6.4]{FreconJaligot07}, it became clear that, for that specific lemma, 
the protocol had implementations autonomous from these specific configurations. 
Proofs may have appeared in \cite{AltinelBurdgesBullShit, Freconpseudotores}, 
with a conceptually better and more general implementation in the second case but, 
regrettably, with no connection at all to \cite{CherlinJaligot2004} in both cases. 

A much more delicate use of the protocol can be found in
\cite[Proposition 6.17]{CherlinJaligot2004}.
It is proved there, still in a specific configuration, that the centralizer of a
certain finite subgroup of a decent torus is connected, with then a much
more restrictive faithful action of the Weyl group. As this special
application of the protocol contains the main difficulty possibly inherent
to the subject, we mostly refer to this example. As we will see below, the
key point is that generosity is in general related to a finiteness
property, as opposed to a uniqueness property, a delicate aspect treated 
``by hand" in \cite[Proposition 6.17]{CherlinJaligot2004} and much more 
conceptually here. 

Theorem \ref{TheoGenerixCosets} has general consequences on the action 
of the Weyl group on the underlying subgroup, again whatever these are. 
Back to the concrete example of a reductive algebraic group, the maximal 
algebraic torus is a divisible abelian subgroup, and the Weyl group acts 
faithfully on it. The main corollary of Theorem \ref{TheoGenerixCosets} is a  
general form of this in the abstract context of groups of finite Morley rank. 

\begin{cor}\label{CorGenWeylGpFaithfull}
Let $G$ be a connected group of finite Morley rank in which, generically, elements 
belong to connected nilpotent subgroups. Suppose that $H$ is a definable connected 
generous subgroup, that $w$ is an element normalizing $H$ but not in $H$, 
of finite order $n$ modulo $H$, 
and that $\{h^{n}~|~h\in H\}$ is generic in $H$. Then $C_{H}(w)<H$. 
\end{cor}

In the case of a connected reductive algebraic group, the subgroup $H$ in
Corollary \ref{CorGenWeylGpFaithfull} is typically the
maximal torus $T$ and $w$ a representative of a nontrivial 
element of order $n$ of the Weyl group. 
In the finite Morley rank case, $H$ may typically be a generous
$n$-divisible Carter subgroup $Q$, and $w$ a representative of a nontrivial element of order 
$n$ of the Weyl group $N(Q)/Q$. One gets then, for instance 
if $Q$ is a divisible abelian generous Carter subgroup as in 
Corollary \ref{CorCarterAbelienGenerousDiv} below, consequences qualitatively 
similar in the finite Morley rank case. 

As for Theorem \ref{TheoGenerixCosets}, the statement adopted
in Corollary \ref{CorGenWeylGpFaithfull} is far more general than what it 
says about this typical case. Less typical applications can be found in \cite[\S4.2]{DeloroJaligotI} 
in the context of connected {\em locally$\o$ solvable$\o$} groups, the smallest class of groups of 
finite Morley rank containing connected solvable groups and Chevalley groups of type $\PSL_{2}$ 
and $\SL_{2}$ over algebraically closed fields. Besides, the reader 
can find there a form of Theorem \ref{TheoGenerixCosets}, actually weaker, but which 
reformats uniformly and in a hopefully informative way the original arguments 
of \cite{CherlinJaligot2004} in this context of ``small" groups. 

\section{Technicalities and environment}

Before passing to the proofs, we review briefly the background needed, or surrounding. 

Groups of finite Morley rank are equipped with a rudimentary
notion of finite dimension on their definable sets, satisfying as axioms
a few basic properties of the natural dimension of varieties in
algebraic geometry
over algebraically closed fields. By definable we mean definable by a
first-order logic formula, possibly with parameters and possibly
in quotients by definable equivalence
relations. The dimension, or ``rank", of a definable set $A$ is denoted
by $\rk(A)$.

The finiteness of the rank implies the descending chain condition on
definable subgroups,
and this naturally gives abstract versions of classical notions of the
theory of algebraic groups:
\begin{itemize}
\item
The {\em definable hull} of an arbitrary subset of the ambient group is
the smallest definable subgroup containing that set. It is contained in
the Zariski
closure in the case of an algebraic group.
\item
The {\em connected component} $G\o$ of a group $G$ of finite Morley rank
is the smallest (normal) definable subgroup of finite index of $G$, and
$G$ is
{\em connected} when $G=G\o$.  
\end{itemize}

A fundamental property of a connected group of finite Morley rank is that it 
cannot be partitioned into two definable {\em generic} subsets, that is two 
subsets of maximal rank \cite{Cherlin79}. Our arguments make full use of the 
following simpler properties. 

\begin{fait}\label{FactGpConActingOnHModH0}\
\begin{itemize}
\item[$(1)$]
A connected group of finite Morley rank acting definably on a finite set must fix 
it pointwise. 
\item[$(2)$]
A connected group of finite Morley rank acting definably on a group $H$
of finite Morley rank induces a trivial action on $H/H\o$. 
\end{itemize}
\end{fait}
\proof
The first item is a well known application of connectedness: as elements of the base set 
have finite orbits, their (definable) stabilizers are of finite index, and hence cannot be proper. 
The second item is a special case of the first which does not seem to be specifically mentioned in the literature: 
as $H\o$ is definably characteristic in $H$, the acting group induces an
action on $H/H\o$, and we are then in presence of the action of a connected group on a
finite set.
\qed

Following \cite{Jaligot06}, we say that a definable subset of a group $G$ of 
finite Morley rank is {\em generous} in $G$ when the union of its $G$-conjugates is 
generic in $G$. In our proof of Theorem \ref{TheoGenerixCosets}, we are essentially
going to reuse lines of arguments of \cite{Jaligot06} for dealing with
generosity, both for characterizing it and for applying it. 
When working with generosity in very general contexts, one has to inspect 
closely each conjugacy class of each individual element of the set considered. The 
reader can find in \cite[\S2.2]{Jaligot06} such an analysis, done there for 
definable connected subgroups. Another approach for this analysis was 
mentioned to the author by Cherlin, with a more conceptual geometric proof, duale in some sense, 
giving also a few more rank equalities. We take here the opportunity to recast these 
computations in terms of permutation groups, not only because it generalizes naturally, but 
also as it certainly might be useful in this more general context. 

Given a permutation group $(G,\Omega)$ and a subset $H$ of $\Omega$,
we denote by $N(H)$ and by $C(H)$ the {\em setwise} and
the {\em pointwise} stabilizer of $H$ respectively, that is
$G_{\{H\}}$ and $G_{(H)}$ in a usual
permutation group theory notation. We also denote by $H^{G}$ the set 
$\{h^{g}~|~(h,g)\in{H \times G}\}$, where $h^{g}$ denotes 
the image of $h$ under the action of $g$, as in the case of an action by conjugation. 
Subsets of the form $H^{g}$ for some $g$ in $G$ are also called
{\em $G$-conjugates} of $H$. Notice that the set $H^{G}$ can be seen, alternatively, as the 
union of $G$-orbits of elements of $H$, or also as the union of $G$-conjugates of $H$. 
When considering the action of a group on itself by conjugation,
as we will do below, all these terminologies and notations are the usual
ones, with $N(H)$ and $C(H)$ the {\em normalizer} and the {\em centralizer} of $H$ respectively.

We note that in this paper we work only with ``exact" normalizers
$N(H)=\{g\in G~|~H^{g}=H\}$, or ``stabilizers", as opposed to ``generic
stabilizers", where the equality $H^{g}=H$ is understood up to a symmetric difference
of lower rank.

\begin{fait}\label{FaitRankComput}
{\bf \cite[Proposition 2.9]{Jaligot06}}
Let $(G,\Omega)$ be a ranked permutation group, $H$ a definable subset
of $\Omega$, and assume that for $r$ between $0$ and $\rk(G/N(H))$ the definable subset
$H_{r}$ of $H$, consisting of those elements of $H$ belonging to a set of $G$-conjugates
of $H$ of rank exactly $r$, is nonempty. Then
$$\rk({H_{r}}^{G})=\rk(G)+\rk(H_{r})-\rk(N(H))-r.$$
\end{fait}
\proof
One may proceed exactly as in the geometric proof of \cite[Proposition
2.9]{Jaligot06}.
In the natural geometry associated to this computation, points are the
elements of
$\Omega$ which are $G$-conjugate to those of $H$ and lines
are the $G$-conjugates of $H$. The set of flags is the set of
couples (point,line) where the point
belongs to the line, and one considers the subflag naturally associated
to $H_{r}$.
Projecting on the set of points one gets $\rk({H_{r}}^{G})+r$ for the rank
of this subflag, and similarly $\rk(G/N(H))+\rk(H_{r})$ by projecting on
the set of lines.
The equality follows.

In this proof we use essentially only two properties of the rank. The first one is a guarantee 
that the sets $H_{r}$ considered are definable. The second one is a 
guarantee of the two formulas as above for the rank of a set, as the sum of the rank of its image 
by a definable function and of the rank of the fibers of that function, when constant. 
These two properties correspond respectively to the {\em definability} and the 
{\em additivity} of the rank in the Borovik-Poizat axioms for ranked structures 
\cite[\S4]{BorovikNesin(Book)94}. 
\qed

\bigskip
In the context of a permutation group as in Fact \ref{FaitRankComput},
we may naturally say that the definable subset $H$ of $\Omega$ is {\em generous} when
the subset $H^{G}$ of $\Omega$ is generic in $\Omega$. Of course, this matches with the usual
definition in the case of the action of a group on itself by conjugation. 
Continuing in the general context of permutation groups, Fact
\ref{FaitRankComput} has the following corollary characterizing generosity.

\begin{cor}\label{CorHGenr=0}
Assume furthermore $\rk(G)=\rk(\Omega)$ and $\rk(H)\leq \rk(N(H))$ 
in Fact \ref{FaitRankComput}. Then $H^{G}$ is generic in $\Omega$ if and only if 
$\rk(H_{0})=\rk(N(H))$. In this case $\rk(H_{0})=\rk(H)=\rk(N(H))$, 
a generic element of $\Omega$ lies in only finitely many conjugates of 
$H$, and the same applies to a generic element of $H$.
\end{cor}
\proof
If $H^{G}$ is generic in $\Omega$, then one has for some $r$ as in
Fact \ref{FaitRankComput} that ${H_{r}}^{G}$ is generic in $\Omega$, and
then
$$0\leq r =\rk(H_{r})-\rk(N(H))\leq \rk(H)-\rk(N(H))\leq 0,$$
showing that all these quantities are equal to $0$. In particular $r=0$, and  
$\rk(H_{0})=\rk(N(H))$. Conversely, if $\rk(H_{0})=\rk(N(H))$,
then $\rk({H_{0}}^{G})=\rk(G)=\rk(\Omega)$ by Fact \ref{FaitRankComput}.

For our last statement, we also see with the above inequalities that
$\rk(H)=\rk(N(H))$, and as $H_{0}$ and $N(H)$ have the same rank it follows
that $\rk(H_{0})=\rk(H)=\rk(N(H))$. In particular the definable subset
$H_{0}$ of $H$
is generic in $H$, and together with the genericity of ${H_{0}}^{G}$ in
$\Omega$ this is
exactly the meaning of our two last claims.
\qed

\bigskip
We stress the fact that, under the circumptances of 
Corollary \ref{CorHGenr=0}, the generosity of $H$ is equivalent to the genericity of 
the {\em definable} sets $H_{0}$ and ${H_{0}}^{G}$ 
in $H$ and $\Omega$ respectively, so that working with these definable
sets avoids troublesome saturation issues. At this point, it is also worth mentioning
that there are uniform bounds on finite sets throughout. 
This is one of the Borovik-Poizat axioms, usually called {\em elimination of infinite quantifiers}, 
which gives uniform bounds on the cardinals of finite sets in uniformly 
definable families of sets. This is used on rather rare occasions, and could also be used 
here to see the definability of sets like $H_{0}$ in Fact \ref{FaitRankComput} and 
Corollary \ref{CorHGenr=0}: $H_0$ is exactly the set of elements of $H$ 
contained in at most $m$ distinct conjugates of $H$, for some {\em fixed} finite $m$. We will 
not use it as the definability of the rank amply suffices here, but this aspect can of 
course be kept in mind. 

A typical case in which Fact \ref{FaitRankComput} and Corollary \ref{CorHGenr=0} apply 
is the case in which the permutation group $(G,\Omega)$ is interpretable in a group $G$ 
of finite Morley rank. In the rest of this paper we are only going to consider the action 
of a group of finite Morley rank on itself by conjugation, so 
Fact \ref{FaitRankComput} and Corollary \ref{CorHGenr=0} will be applied freely.

As $G$ and $\Omega$ are the same is this case, the extra assumption
$\rk(G)=\rk(\Omega)$ is then automatically satisfied in the
characteristion of generosity of Corollary \ref{CorHGenr=0}. 
The second assumption $\rk(H)\leq \rk(N(H))$ is not satisfied in general, but an 
interesting case in which it holds is the case in which $H$ has 
the form $x\Gamma$, where $\Gamma$ is a definable subgroup of $G$ and $x$ is 
an element of $G$ normalizing $\Gamma$: in this case $\Gamma \leq N(x\Gamma)$, and thus
$\rk(x\Gamma)=\rk(\Gamma)\leq \rk(N(x\Gamma))$. In fact, one sees in this case 
that $N(x\Gamma)$ is exactly the preimage in $N(\Gamma)$ of 
$C_{N(\Gamma)/\Gamma}(x\mbox{~mod~}\Gamma)$. All cosets considered in this paper
are of this type, and we will make full use of Corollary
\ref{CorHGenr=0} when considering
the generosity of such cosets in the rest of the paper.

We insist again on the fact that the characterisation of Corollary \ref{CorHGenr=0} is in this 
case essentially the genericity of $H_{0}$ in $H$ (in addition to $\rk(H)=\rk(N(H))$), 
and thus the fact that only finitely many conjugates of $H$ pass through 
a generic element of $H$. In general, and we would like to say with probability 
almost one, there is not uniqueness. It may be seen by considering the generic element $g$ 
of a connected reductive algebraic group. It lies in a maximal torus $T$, which lies in a 
generous Borel subgroup $B$; $T$ is the unique of its conjugates containing $g$ 
(\cite[Corollary 3.8]{Jaligot06}), but there are several conjugates of $B$ 
containing $g$ (and permuted by the Weyl group $N(T)/T$). 

That's all about the background we will use. 
We do not use decent tori and Carter subgroups in the present work, 
Theorem \ref{TheoGenerixCosets} and Corollary \ref{CorGenWeylGpFaithfull}, but, 
as they correspond so closely to its most typical applications, it may be useful to 
recall their definitions and to place more precisely our results in context. 
A decent torus $T$ of a group of finite 
Morley rank is a definable (connected) divisible abelian 
subgroup which coincides with the definable hull of its 
(divisible abelian) torsion subgroup, and a Carter subgroup $Q$ 
is a definable connected nilpotent subgroup of finite index in its 
normalizer (and in particular it satisfies $Q=N\o(Q)$). Both types of subgroups 
exist in any group of finite Morley rank, which is trivial in the first case and follows 
in the second case from a graduated notion of unipotence on certain connected nilpotent 
subgroups, for which decent tori are precisely the first stones \cite[\S3.1]{FreconJaligot07}. 
By \cite{Cherlin05}, maximal decent tori are conjugate in any group of finite Morley 
rank, which indeed follows from the fact that $C\o(T)$ is generous for any 
such decent torus $T$. By \cite{Jaligot06}, generous Carter subgroups are conjugate 
in any group of finite Morley rank. 

We take this opportunity to mention the following corelation between
decent tori and generous Carter subgroups.

\begin{fait}\label{FactToreDescGenCarter}
If $Q$ is a generous Carter subgroup of a group of finite Morley rank, then
$T\leq Q\leq C\o(T)$ for some maximal decent torus $T$, and
$N(T)=C\o(T)\cdot N(Q)$.
\end{fait}

Actually, we will prove something slightly more general than Fact \ref{FactToreDescGenCarter}, 
expanding a bit the existing theory of generous subgroups in passing. 

Recall first that the existence of a generous 
Carter subgroup is, maybe, the main open question at the moment concerning groups of 
finite Morley rank. It is equivalent to the question to know whether any connected group of 
finite Morley rank containing no proper definable connected generous subgroup is nilpotent 
(see \cite[Genericity Conjecture 4.1 b--$\beta$]{Jaligot06}). 
As in \cite[\S4.2]{Jaligot06}, a minimal counterexample to the question of existence of a 
generous Carter subgroup in connected groups has tendency to be {\em semisimple}, i.e., 
with all its normal solvable subgroups trivial, and has no proper definable connected generous 
subgroups. 

\begin{fait}\label{FaitMinimaliteCarters}
Let $G$ be a group of finite Morley rank. 
\begin{itemize}
\item[$(1)$]
If $Q$ is a definable nilpotent subgroup of $G$, then any definable subgroup of $Q$ 
generous in $G$ is of finite index in $Q$. 
\item[$(2)$]
If $Q$ and $H$ are definable subgroups of $G\o$ generous in $G$, with $Q$ nilpotent, then 
$Q\o\leq H\o$ up to conjugacy. 
\item[$(3)$]
If $Q$ is a generous Carter subgroup of $G$, then $Q$ is, up to conjugacy, the unique minimal 
definable subgroup of $G\o$ generous in $G$. 
\end{itemize}
\end{fait}
\proof
$(1)$. Assume $H$ is a definable subgroup of $Q$, generous in $G$. Then $H$ must be of 
finite index in its normalizer, by \cite[Lemma 2.2]{Jaligot06} or more generally 
Corollary \ref{CorHGenr=0}. Now by normalizer condition in infinite nilpotent groups of 
finite Morley rank, $H$ is of finite index in $Q$. 

$(2)$. By Corollary \ref{CorHGenr=0} and connectedness of $G\o$, a generic element 
of $G\o$, say $g$, is in conjugates of $Q$ and $H$, say $Q$ and $H$, and in only finitely many 
such conjugates. Now by \cite[Fundamental Lemma 3.3]{Jaligot06}, 
$N\o(Q\cap H)\leq {N\o(Q)\cap N\o(H)}$, and as 
$N\o(Q)=Q\o$ and $N\o(H)=H\o$ by generosity of $Q$ and $H$ 
(using again \cite[Lemma 2.2]{Jaligot06} or Corollary \ref{CorHGenr=0}), 
we get $N\o(Q\cap H)\leq (Q\cap H)\o$. In particular $Q\cap H$ has finite index 
in its normalizer in $Q$, and is thus of finite index in $Q$ by normalizer condition in 
infinite nilpotent groups of finite Morley rank. In particular, $Q\o\leq H\o$. 

$(3)$. By $(1)$ and connectedness of $Q$, $Q$ is minimal for the generosity of definable 
subgroups of $G\o$. By $(2)$, any definable generous subgroup $H$ of $G\o$ 
contains a conjugate of $Q$, i.e., $Q\leq H\o$ up to conjugacy. 
Hence item $(3)$ follows from the conjugacy of generous Carter subgroups of \cite{Jaligot06}. 
\qed

\bigskip
The core of the proof of Fact \ref{FaitMinimaliteCarters} $(2)$ may seem to be somehow 
hidden in the use of \cite[Fundamental Lemma 3.3]{Jaligot06}, which essentially relies 
on Fact \ref{FactGpConActingOnHModH0} $(1)$. Fortunately, our proof of 
Theorem \ref{TheoGenerixCosets} below will reproduce the content of that lemma, 
with cosets instead of subgroups. 

Fact \ref{FaitMinimaliteCarters} $(3)$ provides a way to see generous Carter subgroups 
in the ostensibly wider class of {\em minimal} definable generous subgroups, where the 
problem of existence somehow shifts to the problem of conjugacy. 

We now add decent tori into the picture. 

\begin{fait}\label{FactArgFrattGen}
Let $G$ be a group of finite Morley rank. 
\begin{itemize}
\item[$(1)$]
If $H$ is a definable generous subgroup of $G\o$, then $H\o$ contains a maximal decent torus $T$ 
of $G$. 
\item[$(2)$]
If $H$ is a definable connected generous subgroup of $G$, minimal with respect to this property, 
and $T$ is a maximal decent torus of $G$ in $H$, then $T\leq Z(H)$.  
\item[$(3)$]
If $T$ is a maximal decent torus and $C\o(T)$ contains a unique minimal definable generous 
subgroup up to conjugacy, say $H$, then $N(T)=C\o(T)\cdot N(H)$. 
\end{itemize}
\end{fait}
\proof
$(1)$. By \cite{Cherlin05}, $C\o(T)$ is generous for any decent torus $T$ of $G$. Arguing as in 
the proof of Fact \ref{FaitMinimaliteCarters} $(2)$, one finds a generic element in 
$C\o(T)\cap H$ and one deduces similarly that $N\o(C\o(T)\cap H)\leq {N\o(C\o(T))\cap N\o(H)}=
{C\o(T)\cap H\o}$. As $T$ is central in $C\o(T)$, this implies in particular that $T\leq H\o$. 

$(2)$. By \cite{Cherlin05}, $C\o_{H}(T)$ is generous in $H$. By transitivity of the generosity 
of definable subgroups \cite[Lemma 3.9 a]{Jaligot06}, one deduces that 
$C\o_{H}(T)$ is generous in $G$, and the minimality of $H$ forces $C\o_{H}(T)=H$, i.e., 
$T\leq Z(H)$. 

$(3)$. We have $T\leq N\o_{C\o(T)}(H)=H\o$ by generosity of $H$ in $C\o(T)$, and thus 
$T\leq Z(H)$. In particular, 
$N(H)\leq N(T)$. Now a Frattini Argument gives the desired decomposition: if $w\in N(T)$, then 
$H$ and $H^{w}$ are two minimal definable generous subgroups of $C\o(T)$, 
$H^{w}=H^{\alpha}$ for some $\alpha$ in $C\o(T)$, and 
$w=w\alpha^{-1}\alpha \in N(H)\cdot C\o(T)$. Notice that $C\o(T)$ is normal in $N(T)$. 
\qed

Fact \ref{FactToreDescGenCarter} follows from 
Facts \ref{FaitMinimaliteCarters} and \ref{FactArgFrattGen}, together with the 
remark that the generous Carter subgroup $Q$ of $G$, containing the maximal decent torus $T$, 
must also be generous in $C\o(T)$ 
(by \cite[Lemma 2.3]{Jaligot06} or Corollary \ref{CorHGenr=0}). 

In presence of a {\em nontrivial} maximal decent torus $T$, the {\em Weyl group} of 
an arbitrary group of finite Morley rank is naturally 
defined as in \cite[Theorem 1.8]{CherlinJaligot2004} as $N(T)/C\o(T)$, 
and in presence of a {\em generous} Carter subgroup $Q$, it is defined as in 
\cite[\S3.3]{Jaligot06} as $N(Q)/Q$. In the first case the original definition 
relied on a particular decent subtorus related to the prime $p=2$, but since the full 
proof of conjugacy of maximal decent tori of \cite{Cherlin05} 
it naturally takes this form. We also mention that the term ``Weyl group" made 
his first appearance, beyond the classical algebraic case, in \cite{Nesin89-a} 
in the context of ``bad" groups of Morley rank $3$, with all possible definitions equivalent 
in this case. 

In Fact \ref{FactToreDescGenCarter}, we 
see that both notions of Weyl group essentially match, with however 
$$N(T)/C\o(T)\simeq (N(Q)/Q)/(N_{C\o(T)}(Q)/Q)$$
isomorphic to a possibly proper quotient of $N(Q)/Q$, and thus a sharper notion 
with the second definition. Hence when both definitions are possible we prefer the second one, 
though the question of equality in general is an interesting issue. 

We note that everything said here with a decent torus $T$ can be stated similarly with 
a pseudo-torus $T$, a slightly more general notion of torus with practically the same properties 
\cite{Freconpseudotores}. 

Besides, we note that \cite{FreconConjCarter07} provides an analysis of non-generous Carter subgroups 
in very specific inductive contexts for groups of finite Morley rank. This yields the conjugacy of 
such non-generous Carter subgroups, and eventually gives in these specific cases the full 
conjugacy of Carter subgroups, in the non-generous case as well as in the generous case. 
In particular, this gives a notion of Weyl group in the most pathological situation in which all 
Carter subgroups would be non-generous, the line antipodal to the one pursed in 
\cite{Jaligot06} and, seemingly, here. 

In Theorem \ref{TheoGenerixCosets} we assume that, {\em generically}, elements of the 
ambient group have a prescribed property: to be in a connected nilpotent subgroup. 
As this property has no first-order character, this can be 
interpreted in two possible ways. It means either that the group is saturated 
and that realizations of the generic type have that property, or, more strongly but with no saturation 
assumption, that the ambient group has a definable generic subset, all of whose elements have 
the property. This ``generic property" is known to be true, in this second form, in the specific case of 
connected {\em locally$\o$ solvable} of finite Morley rank, the smallest 
class containing connected solvable 
groups of finite Morley rank and Chevalley groups of type $\PSL_{2}$ over algebraically 
closed fields (see \cite[Proposition 8.1]{BorovikBurdgesCherlin07}, 
and \cite[\S5.3]{DeloroJaligotI} for an account on this and related topics). 
In any case, the assumption in Theorem \ref{TheoGenerixCosets} is 
much weaker than that of the existence of a generous Carter subgroup, and as the former is known 
in contexts where the latter is not known, it seems relevant at present 
to state Theorem \ref{TheoGenerixCosets}, and its consequences, under this weak assumption. 

\section{Cosets and generosity}

In the present section we pass to the proof of the technical Theorem \ref{TheoGenerixCosets} 
on generous cosets, and in the next we will see its main corollary on Weyl groups. 

In most applications of the general protocol for computing
Weyl groups in groups of finite Morley rank, there is a uniqueness property,
and then rank computations for generosity, or non-generosity, follow more or less immediately 
from the presence of {\em disjoint} unions. We refer for example
to \cite[3.3-3.4]{CherlinJaligot2004}, which was essentially extracted
from the original works on bad groups 
\cite[Theorem 13.3, Claim (d)]{BorovikNesin(Book)94}. 
In general, one can use only finiteness instead of uniqueness for generosity, 
as explained and illustrated abundantly after Corollary \ref{CorHGenr=0}. 
The reader can find in \cite[Proposition 6.17]{CherlinJaligot2004} a concrete application 
of the protocol for Weyl groups which uses finiteness only 
(see actually the preparatory sequence 6.13-6.16, and more specifically 3.16, in that paper), 
and we give here a much more conceptual treatment of this aspect 
via Corollary \ref{CorHGenr=0}. 

Recall that $G$ is a connected group of finite Morley rank in which, generically, elements 
belong to connected nilpotent subgroups, that $H$ is a definable subgroup of $G$ and $w$ is an 
element in $N(H)\setminus H$, and we want to show that $wH$ is not generous in $G$. 

\bigskip
\noindent
{\bf Proof of Theorem \ref{TheoGenerixCosets}.} 
Assume towards a contradiction $wH$ generous in $G$.

We may freely apply Corollary \ref{CorHGenr=0} to the coset $wH$, as remarked 
after that corollary. It follows that $\rk(wH)=\rk(N(wH))$ on the one hand,
and, on the other hand, that $wH$ has a definable generic subset, generous
in $G$, all of whose elements can lie in only finitely many conjugates of $wH$.
In this sense, a generic element $g$ of $G$ is, up to conjugacy, a generic element of
$wH$, and contained in only finitely many conjugates of $wH$.
Of course, $N(wH)\leq N(H)$, and in fact $N(wH)$ is the preimage in
$N(H)$ of $C_{N(H)/H}(w\mbox{~mod~}H)$.
As $H$, $wH$, and $N(wH)$ have the same rank,
$$N\o(wH)=H\o.$$
In particular one sees also that $w$ has finite order modulo $H$.

By assumption, a generic element $g$ of $G$ also belongs to a connected 
nilpotent subgroup $Q$ and, as taking definable hulls does not affect 
connectedness and nilpotence of subgroups in group of finite Morley rank, 
we may assume $Q$ definable. (We note here that the generic property in $G$ holds either 
for the realizations of the generic type in case of saturation of $G$, or on all elements of a 
definable generic subset of $G$, if such a subset exists.) 

Using the connectedness of $G$, one concludes from the two preceding paragraphs 
that a generic element $g$ of $G$ is, on the one hand, in $wH$ (up to conjugacy) and in 
only finitely many of its conjugates, and, on the other hand, in a definable connected 
nilpotent subgroup $Q$. We will get a contradiction from this position of 
tightrope walker of $g$. 

As $g\in wH\cap Q$, we may also assume $w$ in $Q$, replacing the original 
representative $w$ of the coset $wH$ by a representative in $Q$ in necessary. This is 
possible as we may take $g$. Then 
$$wH\cap Q=w(H\cap Q).$$ 
Notice that $w$ still has finite order modulo $H\cap Q$, as the original $w$ had that 
property modulo $H$. The group $\<w\>(H\cap Q)$ is in particular definable, and 
$(H\cap Q)\o$ is exactly its connected component. From now on we concentrate on 
the definable subgroup $\<w\>(H\cap Q)$ of $Q$, and to its normalizer in $Q$. 

$N\o_{Q}(\<w\>(H\cap Q))$ acts by conjugation on the definable subgroup 
$\<w\>(H\cap Q)$. By Fact \ref{FactGpConActingOnHModH0} $(2)$, it induces 
a trivial action on this group modulo its connected component, that is $(H\cap Q)\o$. 
This means that it normalizes each coset of $(H\cap Q)\o$ in $\<w\>(H\cap Q)$. 
In particular, $N\o_{Q}(\<w\>(H\cap Q))$ normalizes the (possibly larger) coset 
$w(H\cap Q)$. 

At this point we use an argument similar to the one used in 
\cite[Fundamental Lemma 3.3]{Jaligot06}. We denote by $X$ the set of elements 
of $w(H\cap Q)$ contained in only finitely many conjugates of $wH$. We note that 
the set $X$ is not empty, as it contains the generic element $g$. We also note that 
the subset $X$ of $wH$ can be contained in only finitely many conjugates 
of $wH$, as it contains the element $g$ which has this property. As 
$N\o_{Q}(\<w\>(H\cap Q))$ normalizes $w(H\cap Q)$, it also normalizes 
$X$, and thus it permutes by conjugation the conjugates of $wH$ 
containing $X$. We are now in presence of the definable action of a connected 
group on a finite set, and it follows from Fact \ref{FactGpConActingOnHModH0} $(1)$ 
that it has a trivial action, or in other words that $N\o_{Q}(\<w\>(H\cap Q))$ 
normalizes each of these finitely many conjugates of $wH$ containing $X$. In particular, 
it normalizes $wH$. 

Hence 
$$N\o_{Q}(\<w\>(H\cap Q)) \leq N\o(wH)=H\o,$$
as noticed earlier, and the definable connected subgroup $N\o_{Q}(\<w\>(H\cap Q))$ 
of $Q$ then satisfies 
$$N\o_{Q}(\<w\>(H\cap Q)) \leq (H\o\cap Q)\o\leq (H\cap Q)\o.$$ 
But as $(H\cap Q)\o$ is exactly the connected component of 
$\<w\>(H\cap Q)$, this inclusion shows that 
$\<w\>(H\cap Q)$ has finite index in its normalizer in $Q$. 
Now definable subgroups of infinite index of 
nilpotent groups of finite Morley rank are of infinite index in their 
normalizers, by the classical finite Morley rank version of the normalizer condition in 
finite nilpotent groups. On finds thus that $\<w\>(H\cap Q)$ has finite index in $Q$, 
and by connectedness of the latter one gets 
$$Q=\<w\>(H\cap Q).$$ 
As $(H\cap Q)$ now has finite index in $Q$, one gets similarly 
$$Q=(H\cap Q).$$ 

At this point one gets a contradiction, either by noticing that $w$ has been pushed 
inside $H$, or that $g$ has been pushed outside $Q$. 
\qed

Theorem \ref{TheoGenerixCosets} has the following slightly more general form, 
where the connectedness 
of the ambient group is dropped and the possibly insinuated saturation assumption 
is slightly weakened. 
We note that in this corollary we do not require the elementary extension to be 
saturated itself, but simply that it is satisfies the same assumption as in 
Theorem \ref{TheoGenerixCosets}. 

\begin{cor}\label{CorGnonConnonSat}
Let $G$ be a group of finite Morley rank having an elementary extension $G^{*}$ in which, 
generically, elements belong to connected nilpotent subgroups. Then the coset $wH$ is not generous 
in $G$ for any definable subgroup $H$ of $G\o$ and any element $w$ in $N_{G\o}(H)\setminus H$. 
\end{cor}
\proof
Assume towards a contradiction $wH$ generous in $G$. As $G$ is a finite 
union of translates of $G\o$, $wH$ is generous in $G\o$. 
As the rank can only go up when passing to an elementary extension, 
one then sees that the canonical extension $[wH]^{*}$ of $wH$, in 
$[G^{*}]\o=[G\o]^{*}$, 
is generous in $[G^{*}]\o$. Now one can apply Theorem \ref{TheoGenerixCosets} in $[G^{*}]\o$. 
\qed

\bigskip
Theorem \ref{TheoGenerixCosets} also has the following desirable application. 

\begin{cor}\label{CorTheoGenHGenHoGen}
Let $G$ be a group of finite Morley rank as in Corollary \ref{CorGnonConnonSat} and let 
$H$ be a definable subgroup of $G\o$. Then $H\setminus H\o$ is not generous 
in $G$ and, if $H$ is generous in $G$, then $H\o$ is generous in $G$, and in fact in any 
definable subgroup containing it. 
\end{cor}
\proof
As $H\setminus H\o$ is a finite union of cosets of $H\o$ normalizing $H\o$, 
the first claim follows from Corollary \ref{CorGnonConnonSat}. Now $H\o$ 
must be generous in $G$ whenever $H$ is, and our last claim is \cite[Lemma 3.9]{Jaligot06} 
or Corollary \ref{CorHGenr=0}. 
\qed

\bigskip
In particular, when Corollary \ref{CorTheoGenHGenHoGen} applies in a connected group 
of finite Morley rank, then the notion of minimal definable generous subgroup, as in 
Facts \ref{FaitMinimaliteCarters} or \ref{FactArgFrattGen}, is the same as 
the notion of minimal definable {\em connected} generous subgroup. 

\section{Cosets and action}

As stressed in the introduction, recovering the action of a Weyl
group on its underlying subgroup from weak information on the elements of the
corresponding cosets is a particularly delicate task. Corollary \ref{CorGenWeylGpFaithfull}
is however a general result of faithfulness following merely from
the nongenerosity provided by Theorem \ref{TheoGenerixCosets}. The
rest of this paper is devoted to the proof of 
Corollary \ref{CorGenWeylGpFaithfull}, or rather of what we see as the 
most interesting intermediary steps.

The most general situation is that of a definable connected generous subgroup $H$, 
and we want to examine the action of $N(H)/H$ on $H$, and much more generally the 
action on $H$ of elements $w$ in $N(H)$. Typically, $H$ may be 
a generous Carter subgroup, with then $N(H)/H$ the natural Weyl group, 
and $w$ a representative of any coset of $H$ in $N(H)$.  

We note that a definable generous subgroup $H$ always satisfies
$$N\o(H)=H\o$$
by Corollary \ref{CorHGenr=0}, and is in particular of finite index in 
its normalizer. We note also that there is a basic result of lifting of torsion in groups of 
finite Morley rank, implying in particular that any element of finite order of 
$N(H)/H$ lifts to an element of $N(H)$ of finite order (and where the primes involved 
in both primary decompositions are the same). In particular, choosing 
an element $w$ of finite order, for example as in Corollary \ref{CorGenWeylGpFaithfull}, 
is always a low cost possibility.  

The following lemma is the natural continuation of
\cite[Lemma 3.4]{CherlinJaligot2004} with the present much better
understanding of generosity as a finiteness property as opposed to a uniqueness property. 
It is the finest corelation one can get between generic elements of
the coset $wH$ and generic elements of $H$ in the typical situation
where the conclusion of Theorem \ref{TheoGenerixCosets} holds. It shows, 
we think, the real power of the method.

\begin{lem}\label{LemmaRelDefHull}
Let $G$ be a group of finite Morley rank, $H$ a definable generous
subgroup of $G$, and $w$ an element in $N(H)\setminus H$ such that
$\<w\>H\setminus H$ is not generous. Then
\begin{itemize}
\item[$(1)$]
The coset $wH$ has a definable subset $[wH]_{\rm{gen}}$, whose
complement is nongeneric in
$wH$, and all of whose elements are in infinitely many conjugates of $wH$.
\item[$(2)$]
The subgroup $H$ has a definable generic subset $H_{\rm{gen}}$ such that, for 
any $x$ in $[wH]_{\rm{gen}}$, the subgroup of $\<w\>H$ containing $x$ and 
defined as 
$$\bigcap_{g\in G,~x\in[wH]^{g}}[\<w\>H]^{g}$$
has an empty intersection with $(H_{\rm{gen}})^G$.
\end{itemize}
\end{lem}

\proof
As $N\o(H)=H\o$ by generosity of $H$ and Corollary \ref{CorHGenr=0}, 
$\rk(wH)=\rk(N(wH))$, and the first claim follows from the nongenerosity 
of $wH$ by Corollary \ref{CorHGenr=0}. 
Again we remark that the sets provided by Corollary \ref{CorHGenr=0} are definable. 

Now one can apply Corollary \ref{CorHGenr=0} to $\<w\>H$ also. The 
generosity of $\<w\>H$ (following that of $H$) then gives a definable 
subset $[\<w\>H]_{0}$, generic in $\<w\>H$, and all of whose elements can lie in 
only finitely many conjugates of $\<w\>H$. If that set had a nongeneric 
intersection with $H$, then it would have a generic intersection with one of the
proper cosets of $H$ in $\<w\>H$, say $w'H$. As all elements lying in this
intersection would be contained in only finitely many conjugates of $w'H$, 
as contained in only finitely many conjugates of $\<w\>H$ and all normalizers 
are finite modulo $H\o$, Corollary \ref{CorHGenr=0} would give the generosity of $w'H$, 
a contradiction to the assumption that $\<w\>H\setminus H$ is not generous. 
One may thus consider a generic element of $H$ as
an element of $H_{\rm{gen}}:=H\cap [\<w\>H]_{0}$, and thus with the
property that it is
in only finitely many conjugates of $\<w\>H$.

Consider now $x$ generic in $wH$ in the sense of the first
claim, i.e., such that $x$ is in infinitely many conjugates of $wH$.
The intersection of subgroups considered in our second claim
is a subgroup of $\<w\>H$. It is contained in infinitely many conjugates $\<w\>H$, 
again as all normalizers are finite modulo $H\o$. 
Hence it contains no conjugates of an element in $H_{\rm{gen}}$, 
as such an element is contained in only finitely many conjugates of $\<w\>H$.
\qed

\bigskip
We mention, parenthetically, that the subgroup as in Lemma \ref{LemmaRelDefHull} $(2)$ 
containing the element $x$ of $wH$ is normalized by $C(x)$. It is definable 
by descending chain condition on definable subgroups, and in particular it contains 
the definable hull of $x$ as a (possibly smaller) subgroup. 

In general, an element $x$ of a coset $wH$ has the form 
$x=wh$ for some $h$ in $H$ and taking powers one gets 
$$(wh)^{n}=w^{n}h^{w^{n-1}}h^{w^{n-2}}\cdots h$$
for any natural number $n$ (some useful formulas when considering torsion 
\cite[\S3.3]{CherlinJaligot2004}). 
Assuming additionally that the element $w$ of $N(H)$ has finite order 
$n$ modulo $H$, which can be done in a general way as explained above, one has 
$$(wh)^{n}=w^{n}h^{n}$$ 
in the easiest case in which $w$ and $h$ commute, with $w^{n}$ in $H$. 
This corelation between the element $wh$ of the coset $wH$ and the $n$-th power 
of the element $h$ of $H$ will be combined to the full force 
of the pure genericity argument of Lemma \ref{LemmaRelDefHull} in our proof of 
Corollary \ref{CorGenWeylGpFaithfull}. 

To this end, our next main step is as follows. 

\begin{lem}\label{LemmeGHw<H}
Let $G$ be a group of finite Morley rank, $H$ a definable connected 
generous subgroup, and $w$ an element in $N(H)$ 
such that $\<w\>H\setminus H$ is not generous in $G$. Then 
$$\{h^{w^{n-1}}h^{w^{n-2}}\cdots h~|~h\in H\}$$ 
is not generic in $H$ for any multiple $n$ of the (necessarily finite) order of $w$ modulo $H$. 
\end{lem}
\proof
Assume towards a contradiction $\{h^{w^{n-1}}h^{w^{n-2}}\cdots h~|~h\in H\}$ 
generic in $H$. Let $\phi~:~wh\mapsto (wh)^{n}$ denotes the definable map, 
from $wH$ to $H$, consisting of taking $n$-powers. As 
$$\phi(wH)=w^{n}\cdot \{h^{w^{n-1}}h^{w^{n-2}}\cdots h~|~h\in H\},$$ 
our contradictory assumption forces that $\phi(wH)$ must be generic in $H$. 

Let $H_{\rm{gen}}$ denote the definable generic subset of $H$ provided by 
Lemma \ref{LemmaRelDefHull} $(2)$. By connectedness of $H$, one gets that 
$H_{\rm{gen}}\cap \phi(wH)$ must be generic in $H$ as well. In particular, 
$\phi^{-1}(H_{\rm{gen}}\cap \phi(wH))$ must be generic in the coset $wH$, and 
one finds an element $x$ in this preimage and in the subset 
$[wH]_{\rm{gen}}$ provided by Lemma \ref{LemmaRelDefHull} $(1)$. 

Now $\phi(x)\in H_{\rm{gen}}$, but as $\phi(x)=x^{n}$, one gets 
$$x^{n}\in H_{\rm{gen}}\cap \<x\>,$$
a contradiction to Lemma \ref{LemmaRelDefHull} $(2)$, as $\<x\>$ is obviously a subgroup 
of the subgroup considered in Lemma \ref{LemmaRelDefHull} $(2)$. 
\qed

\bigskip
Combined with Theorem \ref{TheoGenerixCosets}, one gets the following. 

\begin{cor}\label{Cor13}
Let $G$ be a group of finite Morley rank as in Corollary \ref{CorGnonConnonSat}, 
$H$ a definable connected generous subgroup of $G$, and $w$ an element of $G\o$ in 
$N(H)\setminus H$. Then 
$$\{h^{w^{n-1}}h^{w^{n-2}}\cdots h~|~h\in H\}$$ 
is not generic in $H$ for any multiple $n$ of the (necessarily finite) order of $w$ modulo $H$. 
\end{cor}
\proof
As usual, $H$ is of finite index in its normalizer by Corollary \ref{CorHGenr=0}. 
By Theorem \ref{TheoGenerixCosets}, or rather its slightly more general form, 
Corollary \ref{CorGnonConnonSat}, Lemma \ref{LemmeGHw<H} applies. 
\qed 

\bigskip
If $w$ turned out to centralize $H$ in Lemma \ref{LemmeGHw<H}, then one would get 
$$\{h^{w^{n-1}}h^{w^{n-2}}\cdots h~|~h\in H\}=\{h^{n}~|~h\in H\}$$
and thus Corollary \ref{CorGenWeylGpFaithfull} follows similarly from
Theorem \ref{TheoGenerixCosets} and Lemma \ref{LemmeGHw<H}. Again, 
Corollary \ref{CorGenWeylGpFaithfull} could be stated identically in the slightly more 
general context of groups as in Corollary \ref{CorGnonConnonSat}, taking just care to pick up 
the element $w$ in $G\o$ as in Corollary \ref{Cor13}. 
\qed

\bigskip
Not to come to an abrupt end, we mention the following special case of 
Corollary \ref{CorGenWeylGpFaithfull}, much typical of a connected reductive 
algebraic group, where the maximal torus corresponds to our abelian generous Carter subgroup. 
In this mere application, we do not conclude much more than the faithfulness of the action of the 
Weyl group, but state it in a form emphasizing various subgroups reminiscent of the 
$BN$-pair structure of a reductive algebraic group. 

\begin{cor}\label{CorCarterAbelienGenerousDiv}
Let $G$ be a connected group of finite Morley rank with an abelian generous Carter 
subgroup $Q$, and assume $Q$ $p$-divisible for any prime $p$ dividing the order 
of $N(Q)/Q$. Then $Q$ has (finitely many) proper definable subgroups, corresponding to 
all subgroups of the form $C_{Q}(w)$ for $w$ varying in $N(Q)\setminus Q$, and 
with a canonical definition as the centers of proper cyclic extensions of $Q$ in $N(Q)$. 
In particular, $N(Q)/Q$ acts faithfully on $Q$. 
\end{cor}
\proof
Let $w$ in $N(Q)\setminus Q$, of finite order $n$ modulo $Q$. 
As $Q$ is $p$-divisible for all primes $p$ dividing the order of $N(Q)/Q$, 
its is $n$-divisible, and in particular $Q^{n}=Q$. Now $C_{Q}(w)<Q$ by 
Corollary \ref{CorGenWeylGpFaithfull}. We have shown that 
$C_{Q}(w)<Q$ for any element $w$ in $N(Q)\setminus Q$.  

The fact that there are finitely many possibilities for such subgroups 
$C_{Q}(w)$ follows from their alternative definitions as 
$$C_{Q}(w)=Z(\<w\>Q)$$ 
and from the fact that $N(Q)/Q$ is finite. For a canonical definition of such subgroups, 
one may then take $Z(\<w\>Q)$, with $w$ varying in $N(Q)\setminus Q$. 
\qed

\bigskip
\noindent
{\bf Acknowledgments.} 
El autor aprovecha la oportunidad para agradecer a Nadia M. y su familia 
por una Navidad maravillosa en Galicia durante la cual fue concebido este art\`{\i}culo.

\bibliographystyle{alpha}
\bibliography{biblio}

\end{document}